\documentclass[a4paper,12pt]{article}
\usepackage{hyperref}
\usepackage{amssymb}
\usepackage{amsmath}
\usepackage[latin1]{inputenc} 
\newtheorem{defi}{Definition}[section]
\newtheorem{prop}[defi]{Proposition}
\newtheorem{thm}[defi]{Theorem}
\newtheorem{lem}[defi]{Lemma}
\newtheorem{cor}[defi]{Corollary}
\def\R{{\mathbb R}}
\def\C{{\mathbb C}}

\def\N{{\mathbb N}}
\def\Q{{\mathbb Q}}
\def\F{{\mathbb F}}
\def\OO{{\mathcal O}}

\def\eps{{\varepsilon}}

\newcommand {\tx}[1] {\textrm{#1}}
\newcommand {\resr}[1] {\OO/\pi^{#1}\OO}

\newcommand {\esp}[1]{\underset{ #1}{\mathbb E}}
\newcommand{\qed}{\hfill $\Box$}
\newcommand{\pair}[1]{\langle{#1}\rangle}
\newcommand{\flf}[1]{\lfloor{#1}\rfloor}

\begin{document}

\title{Approximation properties for $p$-adic symplectic groups and lattices}
\author{Benben Liao\footnote{ENS de Lyon and Université de Franche-Comté, supported by ANR OSQPI}} 
\date{\today}
\maketitle
\abstract{
Let $G$ be the symplectic group $Sp_4$ over a non Archimedean local field of {\it any} characteristic. It is proved in this paper that for $p\in[1,4/3)\cup (4,\infty]$ neither the group $G$ nor its lattices have 
the property of approximation by Schur multipliers on Schatten $p$ class
($AP_{pcb}^{Schur}$) of Lafforgue and de~la~Salle.
As a consequence, for any lattice $\Gamma$ in $G,$ 
the associated non-commutative $L^p$ space $L^p(L\Gamma)$ of its von Neumann algebra $L(\Gamma)$ fails 
the operator space approximation property (OAP) and 
completely bounded approximation property 
(CBAP) for $p\in[1,4/3)\cup (4,\infty].$
Together with previous work 
\cite{laff-delasalle,haag-delaat-i,haag-delaat-ii,delaat}, one can conclude that lattices in a higher rank algebraic group over {\it any} local field do not have 
the group approximation property (AP) of Haagerup and Kraus.
It is also shown that 
on some lattice $\Gamma$ in $Sp_4$ over some local field,
the constant function $1$ cannot be approximated by radial functions with bounded (not necessarily completely bounded) Fourier multiplier norms on $C^*_r(\Gamma)$,
nor on $L^p(L\Gamma)$ 
for finite $p>4.$

\section{Introduction}

Let $X$ be a Banach space. Recall that $X$ has the Banach space approximation property (AP), if there exist a net of finite rank operators $T_\alpha\in B(X),$ such that
$\lim_{\alpha}\max_{x\in K}\|T_\alpha x-x\|_B= 0$,
 for any compact subset $K\subset X.$
If furthermore $\sup_\alpha \|T_\alpha\|_{B(X)}<\infty,$ we say that $X$ has bounded approximation property (BAP).
BAP is stronger than AP by definition, and in \cite{groth}
Grothendieck showed that for a reflexive Banach space, AP is equivalent to BAP.

An operator space is a closed linear subspace of bounded linear operators on a Hilbert space. 
An operator space $X\subset B(H)$ is said to have the operator space approximation property (OAP), if there exist a net of finite rank operators $T_\alpha\in B(X),$ such that for any $x\in\mathcal K(\ell^2)\otimes_{min} X\subsetneq B(\ell^2\bar\otimes H),$ we have 
$\lim_\alpha\|Id_{\ell^2}\otimes T_\alpha (x)-x\|_{B(\ell^2\bar\otimes H)} =0.$ 
If moreover, the complete bounded norms of $T_\alpha$ are uniformly bounded
$\sup_\alpha\|T_{\alpha}\|_{cb}<\infty,$ then we say that $X$ has the completely bounded approximation property (CBAP). For an operator space, OAP (resp. CBAP) implies AP (resp. BAP) for the underlying Banach space structure \cite{brown-ozawa}.

Let $\Gamma$ be a countable discrete group. Denote $L(\Gamma)$ its group von Neumann algebra and $L^p(L\Gamma)$ the associated non-commutative $L^p$ space, $p\in[1,\infty).$ 

In \cite{laff-delasalle}, it is shown that for a lattice $\Gamma$ in $SL_3(F)$, where $F$ is any local field (e.g. $\R,\C,\Q_p,\F_p((T))$), $L^p(L\Gamma)$ does not have OAP for $p\in[1,4/3)\cup (4,\infty).$
The result is extended in \cite{delaat} to lattices in 
$Sp(2,\R)$
(i.e. $Sp_4(\R)$)
and $p\in[1,12/11)\cup (12,\infty),$
which is improved in \cite{delaat-delasalle} to $p\in[1,10/9)\cup (10,\infty).$
In this article, we show that $L^p(L\Gamma)$ do not have OAP for lattices $\Gamma$ in $Sp_4(F)$ over {\it any} non Archimedean local field $F,$ and $p\in (1,4/3)\cup (4,\infty).$ 

Following \cite{laff-delasalle}, this is achieved by investigating the property of approximations by Schur multipliers on Schatten class $S^p$ ($AP_{pcb}^{Schur}$) for a locally compact group (see Section~2). 
It enjoys many nice properties: having $AP_{pcb}^{Schur}$ is equivalent for lattices and the ambient group; having $AP_{pcb}^{Schur}$ is equivalent to having $AP_{p'cb}^{Schur}$
where $p'$ is conjugate to $p$; 
for $p=1$ or infinity, it is equivalent to weak amenability.

\begin{thm}\label{ap-p}
Let $F$ be a non Archimedean local field of any characteristic. 
Then neither the symplectic group of 4 by 4 matrices $Sp_4(F)\subsetneq M_{4\times 4}(F)$ nor any of its lattices have the $AP_{pcb}^{Schur}$ for $p\in[1,4/3)\cup (4,\infty].$
\end{thm}

\begin{cor}\label{oap}
Let $F$ be a non Archimedean local field. Then 
for any lattice $\Gamma$ in $Sp_4(F),$
the associated non-commutative $L^p$ space of
its von Neumann algebra $L(\Gamma)$ does not have the operator space approximation property (OAP) nor the completely bounded approximation property (CBAP) for $p\in(1,4/3)\cup (4,\infty).$
\end{cor}
{\bf Remark.} Theorem~\ref{ap-p} and Corollary~\ref{oap} are analogues statements of 
Theorem~D and Theorem~A 
on $SL_3$ in \cite{laff-delasalle}.
From the original results on $SL_3$ one does not see the difference between Archimedean local fields and non Archimedean local fields: the ranges of $p$ obtained in \cite{laff-delasalle} for both cases are the same $(1,4/3)\cup(4,\infty).$
Whereas for $Sp_4,$ the ranges $(1,4/3)\cup (4,\infty)$ obtained in this paper for non Archimedean local fields are better than the ones $(1,10/9)\cup (10,\infty)$ established for Archimedean local fields \cite{delaat-delasalle}.
It is unlikely that this is a genuine difference between local fields, but existed arguments \cite{delaat-delasalle} do not improve the ranges for $Sp_4(\R)$.

As for group approximation properties, recall that
for a discrete group $\Gamma,$ weak amenability for $\Gamma$ is equivalent to CBAP for $C_r^*(\Gamma)$ \cite{haag};
approximation property of Haagerup and Kraus (AP) is equivalent to OAP of $C_r^*(\Gamma)$ 
(Theorem~2.1 \cite{haag-kraus}, see also \cite{brown-ozawa}).

Theorem~\ref{ap-p} together with \cite{laff-delasalle,haag-delaat-i,haag-delaat-ii}, we conclude

\begin{cor}\label{gp-ap}
Let $k$ be a local field, and $G$ be an almost simple algebraic $k$-group with $k$-split rank $\geq 2.$ Then non of the lattices in $G(k)$ has the approximation property (AP) of Haagerup and Kraus \cite{haag-kraus}.
\end{cor}

We turn back to Grothendieck's Banach space AP. 
P.~Enflo constructed the first example of Banach space without AP \cite{enflo}. Later, the natural example of bounded linear operators on a Hilbert space was shown to fail AP \cite{szan}.
It is not known whether there exists countable group $\Gamma$ such that $C^*_r(\Gamma)$ or $L^p(L\Gamma)$ for some finite $p$ fails AP (or even BAP).

Let $\Gamma$ be the finitely generated group $Sp_4(\F_q[T])$ of symplectic matrices over the ring of polynomials $\F_q[T]$ where the coefficients are in the finite field $\F_q$ and $q$ is an odd prime power. It is a lattice in $Sp_4(\F_q((T^{-1}))).$
The following theorem rules out the possibilities of approximations by radial Fourier multipliers on $C^*_r(\Gamma)$ and $L^p(L\Gamma).$

We say that $\ell:\Gamma\to\R_{\geq 0}$ is a length function if $\ell(\gamma\gamma')\leq\ell(\gamma)+\ell(\gamma'),\gamma,\gamma'\in\Gamma.$
A function is called $\ell$-radial if
$f(\gamma)=f(\gamma')$ whenever $\ell(\gamma)=\ell(\gamma')$.

\begin{thm} \label{lattice-Lp}
Let $\Gamma$ be the finitely generated group above.
There exists a length function $\ell:\Gamma\to\R_{\geq 0}$ which is biLipschitz to the word length on $\Gamma,$ such that 
the constant function $1\in C(\Gamma)$ cannot be approximated point-wise by any family of $\ell$-radial 
(not necessarily completely bounded)
Fourier multiplier
$(f_\alpha)_{\alpha\in I}\subset \C\Gamma$ 
on $C^*_r(\Gamma)$ with
  $$\sup_{\alpha\in I}\|m_{f_\alpha}\|_{MC^*_r(\Gamma)}<+\infty,$$
 nor by Fourier multipliers on $L^p(L\Gamma)$ with
 $$\sup_{\alpha\in I}\|m_{f_\alpha}\|_{ML^p(L\Gamma)}<+\infty$$
 for finite $p > 4.$
\end{thm}


As a by-product of the arguments, a similar statement on Schur multipliers on Schatten class is also obtained.
Recall that in \cite{laff-delasalle}, it is shown that for a non discrete group, completely bounded Schur multiplier norms and Schur multiplier norms are equal. Whereas, a conjecture of Pisier postulates that there exists a Schur multiplier on $S^p(\ell^2)$ which is not completely bounded for any finite $p\geq 1.$

\begin{thm}\label{lattice-Sp}
Let $\Gamma$ be the finitely generated group defined above (as in Theorem~\ref{lattice-Lp}).
There exists a length function $\ell$ on $\Gamma$ that is biLipschitz to its word length, such that the following holds: for any $p\in (4,+\infty),$ $1\in C(\Gamma)$ cannot be approximated point-wise by $\ell$-radial functions $f_\alpha\in\C(\Gamma)$ such that their Schur multiplier norms are bounded (not necessarily completely bounded) uniformly
 $$\sup_{\alpha\in I}\|m_{f_\alpha}\|_{MS^p(\ell^2\Gamma)}<+\infty.$$
\end{thm}

The paper is organized as follows.

In Section 2, $AP_{pcb}^{ Schur}$ is recalled, some simple facts about non-commutative $L^p$ spaces and quantitative versions of the theorems above are given (modulo important results in \cite{laff-delasalle}).

In Section 3, the proof of Theorem~\ref{ap-p} is given. The proofs are different for
cases when the characteristic of $F$ is $2$ and when it is different from $2.$
Matrices constructed in \cite{laff-ancienne, sp4-strongt, sp4-obstacle} are used 
and some arguments treating $SL_3$ \cite{laff-delasalle} (in particular Lemma~4.9) can be adapted to the case of $Sp_4$.

In Section 4, the proof of Theorem~\ref{lattice-Lp} is given. 
The reason for restricting to radial functions is technical: the arguments only give estimates for spherical functions on the ambient group. The matrices used in the proof of Theorem~\ref{ap-p} do not apply since they do not give rise to invariant operators. Instead, explicit functions constructed in \cite{laff-orsay,sp4-obstacle} are used in the proof (without using Lemma~4.9 \cite{laff-delasalle}). 

Lastly in Section 5, Theorem~\ref{lattice-Sp} is proved.

{\bf Acknowledgment:} I thank Vincent~Lafforgue for his encouragement to study the problem of group approximation properties for $Sp_4.$
I also thank Mikael~de~la~Salle for numerous helpful discussions and valuable suggestions on several improvements and simplifications of the proofs.

\section{Multipliers on Schatten classes and non commutative $L^p$ spaces}

Let $p\in[1,\infty]$ and $H$ be a Hilbert space. 
For $p<\infty,$
denote $S^p(H)$ the Schatten $p$ class on $H,$ i.e. the subspace of bounded operators $T\in B(H)$ such that the trace $Tr(|T|^p)$ is finite. It is a Banach space with respect to the norm $\|T\|_{S^p(H)}=Tr(|T|^p)^{1/p},T\in S^p(H).$
For $p=\infty,$ denote $S^{\infty}(H)$ the space of compact operators.

Let $X$ be a topological space with a fixed Borel measure. A continuous function $\varphi\in C(X\times X)$ is said to be a Schur multiplier on $S^p(L^2X),$ if for any operator $T\in S^p(L^2X)\cap S^2(L^2X)$ (being a dense subspace of $S^p(L^2X)$) with symbol $(T_{x,y})_{x,y\in X},$ the operator with symbol $(\varphi(x,y)T_{x,y})_{x,y\in X}$ is in $S^p(L^2(X))$
and
$$\|(\varphi(x,y)T_{x,y})_{x,y\in X}\|_{S^p(L^2X)}\leq C\|T\|_{S^p(L^2X)}$$
for some $C>0$ - the smallest $C$ is denoted by $\|\varphi\|_{MS^p(L^2X)}.$
If furthermore there exists some $C'>0$ such that for any operators $(T_{x,y}\in B(H))_{x,y\in X}\in S^p(L^2X\bar \otimes H)$
$$\|(\varphi(x,y)T_{x,y})_{x,y\in X}\|_{S^p(L^2X\bar\otimes H)}\leq 
C'\|T\|_{S^p(L^2X\bar\otimes H)},$$
where $H$ is a Hilbert space,
then we say that $\varphi$ is a completely bounded Schur multiplier on $S^p(L^2X),$ and the smallest possible $C'$ is denoted by $\|\varphi\|_{cbMS^p(L^2X)}.$

Let $G$ be a locally compact group with a fixed Haar measure. 
A continuous function on the group $f\in C(G)$ gives rise to a continuous function $[(x,y)\mapsto f(x^{-1}y)]\in C(G\times G)$
(denoted by $\check f$ in \cite{laff-delasalle})
 on its product $G\times G,$ and if it is a Schur multiplier on $S^p(L^2G)$ then we denote it by $m_f.$ With our notation we have
$$\|m_f\|_{MS^p(L^2G)}=\sup_{T\in S^p(L^2G),\|T\|_{S^p(L^2G)}\leq 1}\|(f(x^{-1}y)T_{x,y})_{x,y\in G}\|_{S^p(L^2G)},$$
and
$$\|m_f\|_{cbMS^p(L^2G)}=\sup\|(f(x^{-1}y)T_{x,y})_{x,y\in G}\|_{S^p(L^2G\bar\otimes H)}$$
where $H$ is a Hilbert space and
the supremum is taken over operators $(T_{x,y}\in B(H))_{x,y\in X}\in S^p(L^2G\bar\otimes H)$ with $\|T\|_{S^p(L^2G\bar\otimes H)}\leq 1.$

\begin{defi}(\cite{laff-delasalle})
Let $G$ be a locally compact topological group with a fixed Haar measure.
Let $p\in[1,\infty]$ as above. 
Say that $G$ has 
the property of approximations by Schur multipliers on Schatten $p$ class
$AP_{pcb}^{Schur},$ if there exist a net of functions $(f_{\alpha})_{\alpha\in I}$ in the Fourier algebra $A(G)$ (being a subset of $C_0(G)$) which are completely bounded Schur multipliers on Schatten $p$ class $S^p(L^2G)$ with uniformly bounded norms
$$\sup_{\alpha \in I}\|m_{f_\alpha}\|_{cbMS^p(L^2G)}< +\infty,$$
such that the constant function $1$ on $G$ can be approximated by these functions $(f_\alpha)_{\alpha\in I}$ uniformly on compact sets.
\end{defi}

Since $S^2(L^2G)$ is the space of Hilbert-Schmidt operators and
$$\|m_f\|_{cbMS^2(L^2G)}=\|f\|_{L^{\infty}(G)},$$
$G$ always has $AP_{2cb}^{Schur}.$

After \cite{bf}, completely bounded multipliers on the Fourier algebra $A(G)$ coincide with that on compact operators on $L^2(G):$
$$\|f\|_{M_0A(G)}=\|m_f\|_{cbMS^\infty(L^2G)},\forall f\in C(G),$$
we see that $AP_{\infty cb}^{Schur}$ (or $AP_{1 cb}^{Schur}$, see \cite{laff-delasalle})
is equivalent to weak amenability for $G.$


Let $\Gamma$ be a countable discrete group. Denote $L(\Gamma)$ its group von Neumann algebra, namely the bicommutant of the operators generated by the left regular representation $\lambda$ of $\Gamma$ on $\ell^2\Gamma$
$$L(\Gamma)=\{\lambda (\gamma),\gamma\in \Gamma\}'' \subset B(\ell^2\Gamma).$$
$L(\Gamma)$ is equiped with the natural faithful tracial state $\tau(x)=\pair{\delta_1,x \delta_1},x\in L(\Gamma).$
For finite $p\geq 1,$
denote $L^p(L\Gamma)$ the non commutative $L^p$ space associated to $L(\Gamma),$ i.e. the Banach space of the completion of $L(\Gamma)$ under the norm
$$\|x\|_{L^p(L\Gamma)}=\big(\tau(|x|^p) \big)^{1/p},x\in L(\Gamma).$$

The following statement is probably well-known to expert.

\begin{prop}\label{Lp-C*}
Let $\Gamma$ be a countable discrete group and $f\in\C(\Gamma).$ We have
$$\|f\|_{L^p(L\Gamma)}\leq \|f\|_{L^q(L\Gamma)},1\leq p\leq q<\infty,$$
and
$$\lim_{p\to\infty}\|f\|_{L^p(L\Gamma)}=\|f\|_{C^*_r(\Gamma)}.$$
\end{prop}
{\bf Proof of Proposition~\ref{Lp-C*}:}

Let $H$ be a Hilbert space and $x\in B(H)$ be a normal operator. Denote $\Omega\subsetneq\C$ the Gelfand spectrum of the abelian $C^*$ algebra generated by $x.$ For any unit vector $\xi\in H,$ by Riesz theorem 
 there exists a Borel probability $\mu$ on $\Omega$ such that 
$$\pair{\xi,F(x)\xi}=\int_\Omega F(t)d\mu,$$
$\forall F\in C(\Omega).$
Now apply it to $x=|f|,F(x)=x^p,\xi=\delta_e\in\ell^2\Gamma,$ and by the inequalities of means we get the results.
\qed




\begin{prop}\label{fg-integral}
Let $H$ be a finite group, $1\leq p<\infty$. Then for any function $f\in\C(H),$ we have
$$|\sum_{h\in H}f(h)|\leq |H|^{1/p}\|f\|_{L^p(LH)}.$$
\end{prop}
{\bf Proof of Proposition~\ref{fg-integral}:}
We first have
$$\|f1_H\|_{L^1(LH)}=\pair{\delta_1,f1_H\delta_1}=\sum_{h\in H}f(h),$$
where $1_H$ denotes the constant function one on $H.$

By Holder inequality 
$$\|f1_H\|_{L^1}\leq \|f\|_{L^p}\|1_F\|_{L^q},$$
where $1/p+1/q=1$ and $\pair{\delta_1,|1_H|^q\delta_1}=|H|^{q-1}$.
\qed

For a countable discrete group $\Gamma$ and $f\in\C\Gamma,$ 
we set $\|f\|_{L^\infty(L\Gamma)}=\|f\|_{C^*_r(\Gamma)},$ and $L^\infty(L\Gamma)=L\Gamma,$
and denote $\|m_f\|_{ML^p(L\Gamma)}$ the Fourier multiplier norm of $f$ on $L^p(L\Gamma),p\in[1,\infty]:$

$$\|m_f\|_{ML^p(L\Gamma)}=\sup_{\varphi\in\C\Gamma,\|\varphi\|_{L^p(L\Gamma)}\leq 1}
\Big\| [\gamma\mapsto f(\gamma)\varphi(\gamma)] \Big\|_{L^p(L\Gamma)}.$$

Now we turn to a quantitative version of Theorem~\ref{ap-p}, based on which the theorem and Corollaries~\ref{oap} and \ref{gp-ap} are direct consequences of results in \cite{laff-delasalle}.

Let $F$ be a non Archimedean local field, $\OO\subset F$ its ring of integer. Let $G=Sp_4(F),$ i.e. the matrices $A\in M_{4\times 4}(F)$ satisfying $A^t JA=J,$ where 
$$J=
\begin{pmatrix}
0&0&0&1
\\0&0&1&0
\\0&-1&0&0
\\-1&0&0&0\end{pmatrix}.$$

\begin{thm} \label{ap-p-decay}
$G=Sp_4(F),K=Sp_4(\OO).$
Let $p\in (4,+\infty].$ There exists a continous function $\phi_p\in C_0(G)$ vanishing at infinity, such that for any $K$-biinvariant continuous function $f\in C(G),$
we have
$$|f(g)|\leq \phi_p(g)\|m_f\|_{MS^p(L^2(G)).}$$
\end{thm}


When $p=\infty,$ the analogue of such an inequality turns out to be important also in the proof of negation of AP for $Sp_4(\R)$ (i.e. $Sp(2,\R)$) in \cite{haag-delaat-i} and property $(T^*)$ for $SL_3(\R)$ and $Sp_4(\R)$ in \cite{haag-delaat-knudby}.

For $p=\infty$ and residue field of $F$ has char different from $2,$ the statement is already known by \cite{laff-ancienne}.

{\bf Proof of Theorem~\ref{ap-p} by Theorem~\ref{ap-p-decay}}:
By Theorem~\ref{ap-p-decay}, $1\in C(G)$ cannot be approximated on compact sets by $K$-biinvariant functions $f_\alpha\in C_0(G)$ with $\sup_{\alpha \in I}\|m_{f_\alpha}\|_{cbMS^p(L^2G)}$ being finite. Since both right and left $K$ actions preserve the norm $\|m_f\|_{cbMS^p(L^2(G))}$ (Proposition~4.2 in \cite{laff-delasalle}), we see that the statement extends to all functions $f_\alpha\in C_0(G)$ and thus $G$ does not have $AP_{pcb}^{Schur}$ for $p\in (4,\infty].$

Now that $AP_{pcb}^{Schur}$ is symmetric for a pair of conjugate numbers $p,p'\in[1,\infty],1/p+1/p'=1$ (Proposition~2.3 in \cite{laff-delasalle}), we have that $G$ fails $AP_{pcb}^{Schur}$ for $p\in[1,4/3)\cup (4,\infty].$

Since $AP_{pcb}^{Schur}$ extends from any lattice to the ambient group (Theorem~2.5 \cite{laff-delasalle}), we conclude that none of the lattices in $G$ has $AP_{pcb}^{Schur}$ for $p\in[1,4/3)\cup (4,\infty].$
\qed

{\bf Proof of Corollaries~\ref{oap} and \ref{gp-ap}:}

Corollary~\ref{oap} is a direct consequence of the following fact: for a discrete group $\Gamma,$ $L^p(L\Gamma)$ having OAP is stronger than $\Gamma$ having $AP_{pcb}^{Schur}$ for the same $p\in(1,\infty)$ (Corollary~3.13 \cite{laff-delasalle}).

The Archimedean case of Corollary~\ref{gp-ap} is proved in \cite{haag-delaat-i,haag-delaat-ii}. For a non Archimedean local field $F$ (and in fact any field), we know that any almost simple algebraic group of split rank $\geq 2$ contains a subgroup that is isomorphic to a quotient of $SL_3(F)$ or $Sp_4(F)$ by a finite normal subgroup \cite{bt,margulis}.
Since for a discrete group, having AP implies having $AP_{pcb}^{Schur}$ for all $p\in (1,\infty)$ (Corollary~3.12 \cite{laff-delasalle}), we conclude the proof by showing the following lemma.

\begin{lem}
Let $G$ be a locally compact group and $N\subset G$ a finite normal subgroup. Let $H$ be the quotient group $H=G/N.$ Let
$f\in C_c(H)\mapsto \tilde f\in C_c(G)$ be the embedding of linear spaces
defined by $\tilde f(g)=f(gN\in H).$ We have
$$\|m_{\tilde f}\|_{cbS^pL^2G}\leq \|m_{f}\|_{cbS^pL^2H}.$$
\end{lem}
Now we prove the lemma.
Let $K$ be a Hilbert space.

Set
$$s^*:B(L^2(G,K))\to B(L^2(H\times N,K))$$
$$(T_{x,y}\in B(K))_{x,y\in G}\mapsto (T_{s(z)n,s(w)m})_{(z,n),(w,m)\in H\times N},$$
where $s:H\to G$ is any fixed section. It is an isometry on the subspace of Schattern class $S^p(BL^2(G,K))$ since $s^*$ is induced from the isomorphism of the underlying Hilbert spaces.

We have $s^*(m_{\tilde f}T)=m_f(s^*(T)),$ since by assumption $N$ is a normal subgroup.
\qed

\qed

The following two theorems are also quantitative versions of Theorem~\ref{lattice-Lp} and Theorem~\ref{lattice-Sp} respectively.

\begin{thm}\label{lattice-Lp-decay}
Let $F=\F_q[[\pi]]$ where $q$ is an odd prime power, and let $\Gamma$ be the lattice $Sp_4(\F_q[\pi^{-1}])$ in $G=Sp_4(F).$
Let $K=Sp_4(\OO).$ 
For any $p\in (4,\infty],$
there exists a function vanishing at infinity $\phi_p\in C_0(\Gamma),$ such that 
for any function $f\in \C(\Gamma)\cap$ $^K C(G)^K$ (i.e. $f\in\C(\Gamma)$ and $f(\gamma)=f(\gamma')$ whenever $\gamma\in \Gamma\cap K\gamma'K, \gamma'\in \Gamma$), we have
$$|f(\gamma)|\leq \phi_p(\gamma)\|m_f\|_{ML^p(L\Gamma)}.$$
\end{thm}

For $p=\infty,$ the statement is a special case of Theorem~1.2 when $s=0$
in \cite{sp4-obstacle}. 

{\bf Proof of Theorem~\ref{lattice-Lp} by Theorem~\ref{lattice-Lp-decay}:}
Denote $D(i,j)$ the diagonal matrix
$$\begin{pmatrix}\pi^{-i} & & & \\
 & \pi^{-j} & & \\
 & & \pi^j & \\
 & & & \pi^i
\end{pmatrix}.$$  

Let $\ell: \Gamma\to\R_{\geq 0}$ be the function defined by $\ell(\gamma)=i$ if $\gamma\in KD(i,j)K,i\geq j\geq 0,$ or equivalently $\ell(\gamma)=\log_q\|\gamma\|$ $=\log_q\max_{1\leq \alpha,\beta\leq 4}|\gamma_{\alpha\beta}|_F.$ It is a length function since 
$(KD(i,j)K)^{-1}=KD(i,j)K$ and $\|g_1\|\|g_2\|\geq\|g_1g_2\|.$ 
It is the length function induced from the Bruhat-Tits building associated to $G,$ and thus
biLipschitz to the word length on $\Gamma$ \cite{lub}.

By definition
$\ell$-radial functions are $K$ biinvariant functions on $\Gamma,$ and this completes the proof.
\qed

\begin{thm}\label{lattice-Sp-decay}
Let $F,\Gamma,G,K$ be as in Theorem.~\ref{lattice-Lp-decay}. Then for any $p\in (4,\infty],$ there exists a function $\phi_p\in C_0(\Gamma)$ such that for any function $f\in\C(\Gamma)\cap$ $^K C(G)^K$ we have
$$|f(\gamma)|\leq \phi_p(\gamma)\|m_f\|_{MS^p(\ell^2\Gamma)}.$$
\end{thm}

{\bf Proof of Theorem~\ref{lattice-Sp} by Theorem~\ref{lattice-Sp-decay}:}
One can take the same length function as in the proof of Theorem~\ref{lattice-Lp} by Theorem~\ref{lattice-Lp-decay}.
\qed

\section{Proof of Theorem~\ref{ap-p-decay}}

Denote $D(i,j)$ the diagonal matrice 
$$\begin{pmatrix}\pi^{-i} & & & \\
 & \pi^{-j} & & \\
 & & \pi^j & \\
 & & & \pi^i
\end{pmatrix}.$$  
The set $\Lambda=\{(i,j)\in\N^2,i\geq j\geq 0\}$ is in bijection with the double cosets $K\backslash G/K$ via 
$(i,j)\mapsto KD(i,j)K.$


\begin{prop}\label{estimate-1} 
Let $G,K,p$ be as in Theorem~\ref{ap-p-decay}.
\begin{itemize}
\item
If the characteristic of $F$ is different from $2,$ denote $v_0\in\N$ the valuation of $2\in F.$
Then we have
for any $K$-biinvariant function $f\in C(G),$
$$|f(D(i,j))-f(D(i,j+1))|\leq 2q^{-\frac{1}{2}(i-j-v_0-1)(1-4/p)}\|m_f\|_{MS^p(L^2(G))},$$
where $(i,j)\in\Lambda$ and $i\geq 1,i-j\geq v_0+1.$
\item
If the characteristic of $F$ is $2,$ then $\forall f\in C(G)$ $K$-biinvariant we have
$$|f(D(i,j))-f(D(i,j+2))|\leq 2q^{-\frac{1}{2}(i-j-2)(1-4/p)}\|m_f\|_{MS^p(L^2(G))},$$
where $i\geq j+2.$
\end{itemize}
\end{prop}

\begin{prop}\label{estimate-2}
Let $F$ be a non Archimdean local field of any characteristic, and $G,K,p$ as in Theorem~\ref{ap-p-decay}.
Let $f$ be any $K$ biinvariant function on $G.$ Then for any $(i,j)\in\Lambda$ with $j\geq 3,$ 
we have
$$|f(D(i,j))-f(D(i+1,j-1))|\leq 2q^2\cdot q^{-(j-2)(1-3/p)}\|m_f\|_{MS^p(L^2(G))}.$$
\end{prop}

{\bf Proof of Theorem~\ref{ap-p-decay} by Propositions~\ref{estimate-1} and \ref{estimate-2} above:} It is similar to the proof of Theorem~3.1 by Proposition~3.2 in \cite{sp4-obstacle}, i.e. a zig-zag argument along the line $i=3j.$
\qed

{\bf Proof of Proposition~\ref{estimate-2}:}

\begin{lem}\label{key-lem}
(lemma 4.9 in \cite{laff-delasalle})
Let $m,n\in\N^*,k\in\{1,2,...,m\}.$ Let $p>2+2/n.$
Let $H$ be a locally compact group, $\alpha,\beta:(\resr{m})^{n+1}\to H$ two maps. 
Let $f\in C_c(H)$ satisfy
$$f(\alpha(a_1,a_2,...,a_n,b)\beta(x_1,x_2,...,x_n,y))=\lambda,$$$$
\tx{ if }y=\sum_{i=1}^n a_ix_i+b+\pi^k\in\resr{m},\tx{ and }$$
$$f(\alpha(a_1,a_2,...,a_n,b)\beta(x_1,x_2,...,a_n,y))=\mu,$$$$
\tx{ if }y=\sum_{i=1}^n a_ix_i+b+\pi^{k-1}\in\resr{m}.$$
Then 
$$|\lambda-\mu|\leq 2q^{-k\eps} \|m_f\|_{MS^p((L^2(H))},$$
where $\eps=\frac{n}{2p}\big(p-(2+\frac{2}{n})\big).$
In particular when $n=1,$
$$|\lambda-\mu|\leq 2q^{-\frac{k}{2}(1-4/p)} \|m_f\|_{MS^p((L^2(H))},$$
when $n=2$
$$|\lambda-\mu|\leq 2q^{-k(1-3/p)} \|m_f\|_{MS^p((L^2(H))}.$$
\end{lem}
\qed

We first prove the case when $char(F)\neq 2.$

We first show that
there exist $k\geq i-j-v_0-1$ and
two maps $\alpha,\beta:(\resr{k})^2\to G$ such that
when $y=ax+b,$ we have 
$$\alpha(a,b)\beta(x,y)\in KD(i,j)K,$$
and when $y=ax+b+\pi^{k-1},$
$$\alpha(a,b)\beta(x,y)\in KD(i,j+1)K.$$
Indeed, one can set  $k=2m-2j-v_0$ where $m$ is the integral part of $(i+j)/2,$
and
$\alpha,\beta$ equal to $\beta^{-1},\alpha$ respectively in the proof of proposition 3.2 in \cite{sp4-strongt}, namely
$$\alpha(a,b)= \begin{pmatrix}
 \pi^{m}&&&\\
&\pi^{i-m+j}&& \\
&&\pi^{-i+m-j}& \\
&&&\pi^{-m}
\end{pmatrix}\cdot\begin{pmatrix}
 1&&&\\
0&1&& \\
\sigma(a)&1&1& \\
\sigma(a)^2-2\sigma(b)&\sigma(a)&0&1
\end{pmatrix},$$ 
$$\beta(x,y)= \begin{pmatrix} 1&&&\\
0&1&& \\
\sigma(x)&0&1& \\
\sigma(x)^2+2\sigma(y)&\sigma(x)&0&1
\end{pmatrix}\cdot\begin{pmatrix} \pi^{-m+j}&&&\\
&\pi^{-m+j}&& \\
&&\pi^{m-j}& \\
&&&\pi^{m-j}
\end{pmatrix},$$ 
where $x,y,a,b\in\resr{k},$ and $\sigma:\resr{k}\to\OO$ is a section.
The computations in \cite{sp4-strongt} show that these matrices indeed satisfy our requirements.
It is also possible to construct $\alpha,\beta$ as variants of the matrices used in \cite{laff-ancienne}.

Now apply Lemma~\ref{key-lem} to $\alpha,\beta$ above, $m=k,$ 
$H=G,$ and $\lambda=f(D(i,j)),\mu=f(D(i,j+1)),$ we have
$$|f(D(i,j))-f(D(i,j+1))|\leq 2q^{-\frac{i-j-v_0-1}{2}(1-4/p)}\|m_f\|_{MS^p(L^2(G))}.$$

Now prove the estimate when $Char(F)=2.$ 

There exist $k\geq (i-j-2)/2$ and $\alpha,\beta:(\resr{k})^2\to G$ such that
when $y=ax+b,$ 
$$\alpha(a,b)\beta(x,y)\in KD(i,j)K,$$
and when $y=ax+b+\pi^{k-1},$
$$\alpha(a,b)\beta(x,y)\in KD(i,j+2)K.$$
We still use the constructions from \cite{sp4-strongt}.
Let $k=m-j-1$ where $m=\flf{\frac{i+j}{2}},$ i.e. the biggest integer $\leq (i+j)/2.$ 
Let $x,y,a,b\in\resr{m-j-1},$ and
$\sigma:\resr{m-j-1}\to\OO$ be a
section, and set
$$\alpha(a,b)= \begin{pmatrix} \pi^{m}&&&\\
0&\pi^{i-m+j}&& \\
\pi^{-i+m-j+1}\sigma(b)&\pi^{-i+m-j}(1+\pi \sigma(a))^2&\pi^{-i+m-j}& \\
0&\pi^{-m+1}\sigma(b)&0&\pi^{-m} \end{pmatrix},$$ 
$$\beta(x,y)=
\begin{pmatrix} \pi^{-m+j}&&&\\
0&\pi^{-m+j}&& \\
\pi^{-m+j}(\sigma(x)+\pi \sigma(y))&0&\pi^{m-j}& \\
\pi^{-m+j}\sigma(x)^2&\pi^{-m+j}(\sigma(x)+\pi \sigma(y))&0&\pi^{m-j}
\end{pmatrix}.$$ 
By similar (or simpler) computations as in the proof of lemme 4.1 in \cite{sp4-strongt} we see that these matrices satisfy our requirements.

Now by applying Lemma~\ref{key-lem} we get
$$|f(D(i,j))-f(D(i,j+2))|\leq 2q^{-\frac{i-j-2}{2}(1-4/p)}\|m_f\|_{MS^p(L^2(G))}.$$
\qed

{\bf Proof of Proposition~\ref{estimate-2}:}
Similarly to the proof of the previous proposition, we will construct appropriate matrices in $G$ and apply Lemma.~\ref{key-lem} to obtain the desired inequality.

When $i+j$ is an even number, there exist $k\geq j-2$ and matrices 
$\alpha,\beta:(\resr{k})^3\to G$ such that $\forall a_1,a_2,b,x_1,x_2,y\in\resr{k},$
if $y=a_1x_1+a_2x_2+b,$ then
$$\alpha(a_1,a_2,b)\beta(x_1,x_2,y)\in KD(i,j)K,$$
and if $y=a_1x_1+a_2x_2+b+\pi^{k-1},$ 
$$\alpha(a_1,a_2,b)\beta(x_1,x_2,y)\in KD(i+1,j-1)K.$$
Indeed, removing the discretization $[\cdot]$ in $\alpha,\beta$ for $i+j\in2\N$ and in $\tilde\alpha,\beta$ for $i+j\in2\N-1$ in 
the proof of the second inequality of proposition 3.2 in \cite{sp4-obstacle} (which are improved constructions of the matrices used in the proof of lemma 2.1 in \cite{laff-ancienne}), we get a construction of $\alpha,\beta.$ More precisely, let $k=m=[(i+j)/2]-1,$ i.e.
when $i+j\in 2\N, m=(i+j)/2-1,$ and when $i+j\in 2\N+1,m=(i+j-1)/2-1.$
 Let $\sigma:\resr{m+1}\to\OO$ be any section. 
 When $i+j$ is even, set
 $\alpha,\beta:(\resr{m+1})^3\to G$ by 
$$\alpha(a_1,a_2,b)=\alpha_1(a_1,a_2,b)$$
$$
=\begin{pmatrix}1&-\pi^{-m-1}(1+\pi\sigma(a_1))&\pi^{-m-1}(1+\pi\sigma(a_2))&-\pi^{-2m}\sigma(b)
\\0&1&0&\pi^{-m-1}(1+\pi\sigma(a_2))
\\0&0&1&\pi^{-m-1}(1+\pi\sigma(a_1))
\\0&0&0&1\end{pmatrix},~~(*)$$
$$\beta(x_1,x_2,y)$$$$
=\begin{pmatrix}1&\pi^{-m}\sigma(x_2)&\pi^{-m}\sigma(x_1)&
\pi^{-m-1}\pi^{-m}(\sigma(x_1)+\sigma(x_2))+\pi^{-2m}\sigma(y)
\\0&1&0&\pi^{-m}\sigma(x_1)
\\0&0&1&-\pi^{-m}\sigma(x_2)
\\0&0&0&1\end{pmatrix}.$$
And when $i+j$ is odd, set 
$$\alpha(a_1,a_2,b)=
\begin{pmatrix}
1&0&0&0\\
&1&\pi^{-1}&0\\
&&1&0\\
&&&1
\end{pmatrix}
\alpha_1(a_1,a_2,b),$$
where $\alpha_1$ is as defined in $(*).$
Identical (after removing $[\cdot]$) computations as in \cite{sp4-obstacle} show that they satisfy required properties. Note that even though in \cite{sp4-obstacle} the local field $F$ is assumed to have characteristic different from $2,$ the constructions of $\alpha,\tilde\alpha,\beta$ are valid for  any characteristic.

Now apply Lemma~\ref{key-lem} to $k,\alpha,\beta,$ and $H=G,n=2, \lambda=f(D(i,j)),\mu=f(D(i+1,j-1)),$ 
we have 
$$|f(D(i,j))-f(D(i+1,j-1))|\leq 2q^{-(j-2)(1-3/p)}\|m_f\|_{MS^p(L^2(G))}.$$ 
\qed

\section{Proof of Theorem~\ref{lattice-Lp-decay}}

We adopt the notations $F,\OO,G,K,D(i,j),\Lambda$ as in Section~3.
Note that the ring of integer $\OO$ is $\F_q[[\pi]].$

\begin{prop}\label{estimate-lattice-Lp}
Let $F,G,K,\Gamma,p$ be as in Theorem~\ref{lattice-Lp-decay}. Then for any function $f\in^KC(G)^K$ we have
$$|f(D(i,j))-f(D(i,j+1))|\leq C_{q,p}q^{-(1/2-2/p)(i-j)}\|m_f\|_{ML^p(L\Gamma)},$$
and 
$$|f(D(i,j))-f(D(i+1,j-1))|\leq C_{q,p} q^{2(i+j)/p-j}\|m_f\|_{ML^p(L \Gamma)}.$$
\end{prop}

{\bf Proof of Theorem~\ref{lattice-Lp-decay} using Proposition~\ref{estimate-lattice-Lp}:}
For any $p>4,$ there exists $n\in\N$ such that $2(1+1/n+1)/p-1<0.$
A zig-zag argument near the line $i=(1+1/n)j$ will yield the estimate.
\qed

{\bf Proof of Proposition~\ref{estimate-lattice-Lp}:}

\begin{lem}\label{explicit}
For each $(i, j)\in\Lambda,$
there exist two finite subgroups $H_{1,i,j},H_{2,i,j}\subsetneq\Gamma$ of cardinality $q^{2(i-j)+3}$ and $q^{2(i+j)+2}$ respectively, and
two family of functions $h_{1,i,j},h_{1,i,j+1}\in\C H_{1,i,j},$ $h_{2,i,j},h_{2,i+1,j-1}\in\C H_{2,i,j}$ 
that are normalized characteristic functions of points in
$KD(i,j)K\cap H_{1,i,j}$ and $KD(i,j)K\cap H_{2,i,j}$ respectively, such that
$$\|h_{1,i,j}-h_{1,i,j+1}\|_{C^*_r(H_{1,i,j})}\leq 2q^{-(i-j)/2},$$
and
$$\|h_{2,i,j}-h_{2,i+1,j-1}\|_{C^*_r(H_{2,i,j})}\leq 2q^2q^{-j}.$$
\end{lem}
Now prove the first inequality.
We set

$$H_{1,i,j}=\{\alpha(a,b,\eps)=\begin{pmatrix}1 &0&\pi^{-i}a&\pi^{-i}b\\
 & 1 &\pi^{-i}\eps&\pi^{-i}a\\
 & & 1 &0\\
 & & & 1
\end{pmatrix},a,b\in\F_q+\F_q\pi+...+\F_q\pi^{i-j}, \eps\in\F_q\},$$
and the following function 
$$h_{1,i,j}=\esp{a\in\F_q+...+\F_q\pi^{i-j}}e_{\alpha(a,t_{i-j}(a^2)+\pi^{i-j},1)},$$
where $t_{i-j}:\F_q+...+\F_q\pi^{2i-2j}\to\F_q+...+\F_q\pi^{i-j}$ is the obvious truncation
$t_{i-j}(\sum_{k\geq 0}a_k\pi^k)=\sum_{0\leq k\leq i-j}a_k\pi^k.$

Let $\chi\in\hat H_{1,i,j},$ and suppose $\chi_1,\chi_2$ are characters of $\F_q+...+\F_q\pi^{i-j}$ and $\chi_3\in\hat{\F}_q$ such that $\chi(\alpha(a,b,\eps))=\chi_1(a)\chi_2(b)\chi_3(\eps).$ 
We have the following:
if $\chi(h_{1,i,j}-h_{1,i,j+1})\neq 0$, then there exists $\theta\in\F_q+...+\F_q\pi^{i-j}$ such that $\chi_1(a)=\chi_2(t_{i-j}(\theta a)).$
Indeed, if $k_\alpha,\alpha=1,2$ is the smallest integer $k$ such that $\chi_\alpha$ is trivial on $\F_q\pi^{i-j-k}+\F_q\pi^{i-j-k+1}+...+\F_q\pi^{i-j}$ and non-trivial on $\F_q\pi^{i-j-k-1}$, then we have $k_1\geq k_2$ unless $\chi(h_{1,i,j}-h_{1,i,j+1})=0.$ The existence of $\theta$ follows from the fact that $\F_q[\pi]/\pi^{i-j+1}\F_q[\pi]$ is a local ring.

By a lemma on Gauss sum \cite{laff-orsay}(see also Lemma~4.3 \cite{sp4-obstacle}) we have that 
$$|\chi(h_{1,i,j}-h_{1,i,j+1})|\leq 2q^{-(i-j)/2},\forall \chi\in \hat{H}_{1,i,j}.$$
This yields the first inequality.
\footnote{
If we set 
$$h'_{1,i,j}=\esp{a,b,c\in\resr{i}}
e_{h_1([\pi^{-i}a],\pi^{-i},[\pi^{-i}a^2+\pi^{-j}(1+\pi c)])},$$
then we also have $\|\Delta=h'_{1,i,j}-h'_{1,i,j+1}\|_{C^*_r (H'_{1,i,j})}\leq 2q^{-(i-j)/2}$ for some finite abelian subgroup $H'_{1,i,j}$ and the support of the spectrum of $\Delta$ has cardinality $\leq q^{2i-2j}.$ Since what contributes in the $L^p$ norm of the spectrum of $\Delta$ is the measure of its support, this gives a second proof of the first estimate in Proposition~\ref{estimate-lattice-Lp}.
}

For the second inequality in the lemma, set
$$H_{2,i,j}=\{\beta(a,b,c)=\begin{pmatrix}1 &[\pi^{-m}a]&[\pi^{-m}b]&[\pi^{-2m}c]\\
 & 1 &0&[\pi^{-m}b]\\
 & & 1 &-[\pi^{-m}a]\\
 & & & 1
\end{pmatrix},a,b,c\in\OO\},$$
where $m=[(i+j)/2],$ and $[\cdot]:\F_q((\pi))\to\F_q[\pi^{-1}]$ is defined by taking the integral part 
$[\sum_{i}a_i\pi^{-i}]=\sum_{i\geq 0}a_i\pi^{-i}.$

The constructions of $h_{2,i,j}$ are identical to the explicit functions $h_{2,i,j}$ used in the proof of Proposition~4.1 in \cite{sp4-obstacle}, namely
$$h_{2,i,j}=\esp{a,b,c\in\resr{i}}e_{\beta(1+\pi a,b/2, \pi^{2m-i}(1+\pi c))}.$$
The second inequality is exactly the second inequality of Proposition~4.2 \cite{sp4-obstacle}.
\qed

Let us first show the first estimate. Apply Proposition~\ref{fg-integral} to $H=H_{1,i,j}$ and $\phi=m_f(h_{1,i,j}-h_{1,i,j+1}),$ 
and by Lemma~\ref{explicit} $h_{1,i,j}$ is a normalized characteristic function supported on $KD(i,j)K,$
we have
$$|f(D(i,j))-f(D(i,j+1))|\leq q^{(2(i-j)+3)/p}\|m_f(h_{1,i,j}-h_{1,i,j+1})\|_{L^p(L(H_{1,i,j}))},$$
and again by Lemma~\ref{explicit} and Proposition~\ref{Lp-C*} it is
$$\leq q^{(2(i-j)+3)/p}\|m_f\|_{ML^p(L\Gamma)}\|h_{1,i,j}-h_{1,i,j+1}\|_{L^p(L(H_{1,i,j}))}$$
$$\leq C_{q,p}q^{(2i-2j)/p}\|m_f\|_{ML^p(L\Gamma)}\|h_{1,i,j}-h_{1,i,j+1}\|_{C^*_r(H_{1,i,j})}$$
$$\leq C_{q,p}\|m_f\|_{ML^p(L\Gamma)}q^{-(i-j)(1/2-2/p)}.$$

The second estimate is proved in the same way.
\qed

\section{Proof of Theorem~\ref{lattice-Sp-decay}}

We adopt the notations $F,\OO,G,K,D(i,j),\Lambda$ as in Section~3.

\begin{prop}\label{estimate-lattice-Sp}
Let $F,G,K,\Gamma,p$ be as in Theorem~\ref{lattice-Sp-decay}. Then we have for any function $f\in\C(\Gamma)\cap {^KC(G)^K}$
$$|f(D(i,j))-f(D(i,j+1))|\leq 2q^{-\frac{1}{2}(i-j-2)(1-4/p)}\|m_f\|_{MS^p(\ell^2 \Gamma)},$$
and 
$$|f(D(i,j))-f(D(i+1,j-1))|\leq 2 q^{-(j-2)(1-3/p)}\|m_f\|_{MS^p(\ell^2 \Gamma)}.$$
\end{prop}
We remark that the arguments in Section~4 yield the same (up to a constant) decaying factor $q^{(i-j)(1/2-2/p)}$ for the first inequality and a worse one $q^{2(i+j)/p-j}(>q^{-j(1-3/p)})$ for the second inequality. To be consistant a complete proof of the first inequality is also given below.

{\bf Proof of Proposition~\ref{estimate-lattice-Sp}:}
The proof proceeds in a similar way as the proof of Propositions~\ref{estimate-1} and \ref{estimate-2} - we will construct matrices satisfying required conditions and the apply Lemma~\ref{key-lem}.

We prove prove the first estimate.

There exist $\alpha,\beta:(\resr{i+1})^2\to \Gamma$ such that for $y=ax+b$
we have
$$\alpha(a,b)\beta(x,y)\in KD(i,j)K\cap\Gamma$$
and for 
$y=ax+b+\pi^{i-j-1}$
$$\alpha(a,b)\beta(x,y)\in KD(i,j+1)K\cap\Gamma.$$
The construction of $\alpha,\beta$ are identical to $\alpha,\beta$ used in the first proof of Proposition~3.2 in \cite{sp4-obstacle} (which are matrices in \cite{laff-ancienne} after discretization). More precisely, set
$$\alpha(a,b)
=\begin{pmatrix}1&0&[\pi^{-i}\sigma(a)]&[\pi^{-i}\sigma(a^2-b)]
\\0&1&\pi^{-i}&[\pi^{-i}\sigma(a)]
\\0&0&1&0
\\0&0&0&1\end{pmatrix},$$
$$\beta(x,y)
=\begin{pmatrix}1&0&[\pi^{-i}\sigma(x/2)]&[\pi^{-i}\sigma(x^2/4+y)]
\\0&1&0&[\pi^{-i}\sigma(x/2)]
\\0&0&1&0
\\0&0&0&1\end{pmatrix},a,b,x,y\in\resr{i+1}$$
where $[\cdot]:\F_q((\pi))\to\F_q[\pi^{-1}]$ the integral part of an element
and
$\sigma:\resr{i+1}\to\OO$ is a section. 

Apply Lemma~\ref{key-lem} to $H=\Gamma,\alpha,\beta,k=i-j,m=i+1,\lambda=f(D(i,j)),\mu=f(D(i,j+1))$ we get
$$|f(D(i,j))-f(D(i,j+1))|\leq 2q^{-\frac{1}{2}(i-j-2)(1-4/p)}\|m_f\|_{MS^p(\ell^2 \Gamma)}.$$

Now prove the second inequality.

There exist $k\geq j-2$ and 
$\alpha,\beta:(\resr{k})^3\to G$ 
such that when $y=a_1x_1+a_2x_2+b,\forall a_1,a_2,b,x_1,x_2,y\in\resr{k}$ we have
$$\alpha(a_1,a_2,b)\beta(x_1,x_2,y)\in KD(i,j)K,$$
and when $y=a_1x_1+a_2x_2+b+\pi^{k-1},$ 
$$\alpha(a_1,a_2,b)\beta(x_1,x_2,y)\in KD(i+1,j-1)K.$$
The constructions of $\alpha,\beta$ are identical to $\alpha,\beta$ in the proof of the second inequality of proposition 3.2 in \cite{sp4-obstacle} when $i+j$ is an even number, and identical to $\tilde\alpha,\beta$ when $i+j$ is odd. They are already used in the proof of Proposition~\ref{estimate-2} which we omit here.

By applying Lemma~\ref{key-lem} to $H=\Gamma,n=2,\alpha,\beta,k,\lambda=f(D(i,j)),\mu=f(D(i+1,j-1))$ we have
$$|f(D(i,j))-f(D(i+1,j-1))|\leq 2 q^{-(j-2)(1-3/p)}\|m_f\|_{MS^p(\ell^2 \Gamma)}.$$
\qed


\begin{thebibliography}{GdlH90}

\bibitem[BT]{bt}
A.~Borel and J.~Tits. 
\newblock Groupes Réductifs. 
\newblock {\em Publ. Math. l'I.H.E.S.} 27 (1965) 55-151.

\bibitem[BF]{bf}
M.~Bozejko and G.~Fendler. 
\newblock Herz-Schur multipliers and completely bounded multipliers of the Fourier algebra of a locally compact group, 
\newblock {\em Boll. Un. Mat. Ital.} A (6) 3 (1984), 297-302. MR 0753889.

\bibitem[BO]{brown-ozawa}
N. P. Brown and N. Ozawa.
\newblock $C^*$-Algebras and Finite-Dimensional Approximations, Grad. Stud. Math. 88.
\newblock {\em Amer. Math. Soc.}, Providence, 2008. MR 2391387.


\bibitem[Enf]{enflo}
P.~Enflo. 
\newblock A counterexample to the approximation property in Banach spaces. 
\newblock {\em Acta Math.,} 130 (1973), 309-317.

\bibitem[Gro]{groth}
A.~Grothendieck. 
\newblock Produits tensoriels topologiques et Espaces nucléaires.
\newblock {\em Mem. Amer. Math. Soc.} 1955 (1955), no. 16, 140 pp. 

\bibitem[Haag]{haag}
U.~Haagerup. 
\newblock Group $C^*$-algebras without the completely bounded approximation property. 
\newblock Preprint, 1986.

\bibitem[HdL13a]{haag-delaat-i}
U.~Haagerup and T.~de~Laat.
\newblock Simple Lie groups without the approximation property.
\newblock {\em Duke Math. J.} Volume 162, Number 5 (2013), 925-964.

\bibitem[HdL13b]{haag-delaat-ii}
U.~Haagerup and T.~de~Laat.
\newblock Simple Lie groups without the approximation property II.
\newblock Trans. Amer. Math. Soc.

\bibitem[HKdL]{haag-delaat-knudby}
U.~Haagerup, S.~Knudby and T.~de~Laat.
\newblock A complete characterization of connected Lie groups with the Approximation Property.
\newblock http://arxiv.org/abs/1412.3033

\bibitem[HK]{haag-kraus}
U.~Haagerup and J.~Kraus.
\newblock Approximation properties for group $C^*$-algebras and group von Neumann algebras.
\newblock {\em Trans. Amer. Math. Soc.} 344 (1994), 667-699. MR 1220905

\bibitem[LMR]{lub}
A.~Lubotzky, S.~Mozes and M.S.~Raghunathan.
\newblock 
The word and Riemannian metrics on lattices of semisimple groups.
\newblock
{\em Publications Mathématiques de l'Institut des Hautes Études Scientifiques},
December 2000, Volume 91, Issue 1, pp 5-53.

\bibitem[dL]{delaat}
T.~de~Laat.
\newblock Approximation properties for noncommutative $L_p$-spaces associated with lattices in Lie groups.
\newblock {\em J. Funct. Anal.} 264 (2013), no. 10, 2300-2322.

\bibitem[dLdlS]{delaat-delasalle}
T.~de~Laat and M.~de~la~Salle.
\newblock Strong property (T) for higher rank simple Lie groups.
\newblock http://arxiv.org/abs/1403.6415

\bibitem[Laf10a]{laff-orsay}
V.~Lafforgue.
\newblock Strong property (T) and the Baum-Connes conjecture, unpublished notes.
\newblock École thématique autour de la conjecture de Baum-Connes et du principe d'Oka en géométrie non-commutative, 2010, Département de Mathématiques d'Orsay, Université de Paris-Sud 11.


\bibitem[Laf10b]{laff-ncg}
V.~Lafforgue.
\newblock Propriété (T) renforcée et conjecture de Baum-Connes, Quanta of maths, 323-345, 
\newblock Clay Math. Proc., 11, Amer. Math. Soc., Providence, RI, 2010

\bibitem[LdlS]{laff-delasalle}
V.~Lafforgue, M.~De~la~Salle.
\newblock Noncommutative $L^p$-spaces without the completely bounded approximation property
\newblock {\em Duke Math. J.} Volume 160, Number 1 (2011), 71-116.

\bibitem[Laf10c]{laff-ancienne}
V.~Lafforgue.
\newblock Un analogue non archimédien d'un résultat de Haagerup et lien avec la propriété (T) renforcée (ancienne version de 2010).
\newblock http://vlafforg.perso.math.cnrs.fr/haagerup-rem-2010.pdf

\bibitem[Liao13]{sp4-strongt}
B.~Liao.
\newblock Strong Banach property (T) for simple algebraic groups of higher rank.
\newblock {\em J. Topol. Anal.}, 06, 75 (2014). DOI: 10.1142/S1793525314500010

\bibitem[Liao14]{sp4-obstacle}
B.~Liao.
\newblock About the difficulty to prove the Baum Connes conjecture without coefficient for a non-cocompact lattice in $Sp_4$ in a local field.
\newblock Submitted. http://arxiv.org/abs/1411.6151


\bibitem[Laf08]{laff-duke}
V.~Lafforgue.
\newblock
 Un renforcement de la propriété $(T)$.
\newblock {\em Duke Math. J.}, 143(3):559--602, 2008.

\bibitem[Laf09]{laff-jta}
V.~Lafforgue.
\newblock
 Propriété $(T)$ renforcée banachique et transformation de Fourier rapide.
\newblock {\em Journal of Topology and Analysis}, Volume: 1, Issue: 3(2009) pp. 191-206.

\bibitem[Mar]{margulis}
G.~A.~Margulis.
\newblock Discrete Subgroups of Semisimple Lie Groups. 
\newblock Springer-Verlag, 1991.


\bibitem[Szan]{szan}
A.~Szankowski.
\newblock $B(H)$ does nol have the approximation properly. 
\newblock {\em Acta Math.}147 (1981),89-108.



\end{thebibliography}
\end{document}